\documentclass[11pt]{article}
\usepackage{amsmath,amsthm}
\usepackage{amssymb,mathrsfs}
\usepackage{bm,mathtools}
\usepackage{tikz}
\usetikzlibrary{calc}

\usepackage{xcolor}
\usepackage{geometry}
\usepackage[colorlinks=true,
linkcolor=blue,citecolor=blue,
urlcolor=blue]{hyperref}

\makeatletter
\def\@seccntDot{.}
\def\@seccntformat#1{\csname the#1\endcsname\@seccntDot\hskip 0.5em}
\renewcommand\section{\@startsection{section}{1}{\z@}%
{18\p@ \@plus 6\p@ \@minus 3\p@}%
{9\p@ \@plus 6\p@ \@minus 3\p@}%
{\large\bfseries\boldmath}}
\renewcommand\subsection{\@startsection{subsection}{2}{\z@}%
{12\p@ \@plus 6\p@ \@minus 3\p@}%
{3\p@ \@plus 6\p@ \@minus 3\p@}%
{\bfseries\boldmath}}
\renewcommand\subsubsection{\@startsection{subsubsection}{3}{\z@}%
{12\p@ \@plus 6\p@ \@minus 3\p@}%
{\p@}%
{\bfseries\boldmath}}
\makeatother

\usepackage{microtype}

\theoremstyle{plain}
\newtheorem{theorem}{Theorem}

\newtheorem{conjecture}{Conjecture}
\newtheorem{problem}{Problem}

\theoremstyle{definition}

\newtheorem{remark}{Remark}

\allowdisplaybreaks
\DeclareMathOperator{\spex}{spex}
\DeclareMathOperator{\SPEX}{SPEX}

\title{Unsolved Problems in Spectral Graph Theory}

\author{Lele Liu\footnote{School of Mathematical Sciences, Anhui University, Hefei 230601, 
P.R. China. E-mail: \texttt{liu@ahu.edu.cn} (L. Liu). Supported by the National
Nature Science Foundation of China (No. 12001370)}
~~ and ~~
Bo Ning\footnote{Corresponding author. College of Computer Science, Nankai University, Tianjin 300350, P.R. China.
E-mail: \texttt{bo.ning@nankai.edu.cn} (B. Ning). Partially supported by the National Nature Science
Foundation of China (No. 11971346).}}
\date{}

\begin{document}
\maketitle

\begin{abstract}
Spectral graph theory is a captivating area of graph theory that employs the eigenvalues 
and eigenvectors of matrices associated with graphs to study them. In this paper, we 
present a collection of $20$ topics in spectral graph theory, covering a range of open 
problems and conjectures. Our focus is primarily on the adjacency matrix of graphs, 
and for each topic, we provide a brief historical overview.
\end{abstract}


\vspace{10mm}

Spectral graph theory is a beautiful branch of graph theory that utilizes the eigenvalues 
and eigenvectors of matrices naturally associated with graphs to study them. The primary 
objective in spectral graph theory is twofold: firstly, to compute or estimate the eigenvalues 
of these matrices and secondly, to establish links between the eigenvalues and the 
structural properties of graphs. As it turns out, the spectral perspective is a powerful 
tool in the study of graph theory.

Over the past few decades, numerous results and applications in various fields of mathematics
have been obtained through spectral graph theory. However, many open problems and conjectures 
in spectral graph theory remain unresolved, necessitating further exploration.

In this paper, we collect $20$ topics in spectral graph theory that include a range of 
conjectures and open problems, with a focus primarily on the adjacency matrix of graphs.
Additionally, we provide a brief historical overview of each topic. Inevitably, it is a somewhat 
personal perspective on the choice of these problems and conjectures, and is not intended 
to be exhaustive; the authors apologize for any omissions.

Let us begin with some definitions and notation. Throughout this paper, we only consider 
graphs that are simple (there are no loops or multiple edges), undirected and unweighted. 
Given a graph $G$, the \emph{adjacency matrix} $A(G)$ of $G$ is a $(0,1)$-matrix, where the 
rows and columns are indexed by the vertices in $V(G)$. The $(i,j)$-entry of $A(G)$ is 
equal to $1$ if the vertices $i$ and $j$ are adjacent, and $0$ otherwise. Since $A(G)$ is 
a real and symmetric matrix, it has a full set of real eigenvalues which we will denote by
\[
\lambda_1(G)\geq \lambda_2(G)\geq\cdots \geq \lambda_n(G).
\]
Recall that the Laplacian matrix of $G$ is defined as $L(G) := D(G) - A(G)$, where $D(G)$ is 
the diagonal matrix whose entries are the degrees of the vertices of $G$. We shall write
$\mu_1(G)\geq\mu_2(G)\geq\cdots\geq\mu_n(G)=0$ for the eigenvalues of $L(G)$. If there is no
danger of ambiguity, for any $1\leq i \leq n$ we write $\lambda_i$ and $\mu_i$ instead of
$\lambda_i(G)$ and $\mu_i(G)$, respectively. We also write $\lambda(G):=\lambda_1(G)$ for short.

For a graph $G$, let $\omega(G)$ denote the clique number of $G$, which is defined as the number
of vertices in the largest complete subgraph in $G$. Let $e(G)$ and $\overline{G}$ denote the 
number of edges and the complement of $G$, respectively. We use the notations $\delta(G)$, 
$\Delta(G)$, and $\overline{d}(G)$ to represent, respectively, the minimum degree, maximum 
degree, and average degree of $G$. For two vertex-disjoint graphs $G$ and $H$, we 
use $G \vee H$ to denote the join of $G$ and $H$, which is obtained by adding all possible 
edges between $G$ and $H$. The \emph{complete split graph} $S_{n,k}$ with parameters $n$ 
and $k$ is the graph on $n$ vertices obtained from a clique on $k$ vertices and an independent 
set on the remaining $n - k$ vertices in which each vertex of the clique is adjacent to each 
vertex of the independent set. If there is no special explanation, we use $n$ to denote 
the number of vertices in $G$, and $m$ to denote the number of edges in $G$. As usual,
$K_n$, $K_{p,n-p}$, $P_n$ and $C_n$ denote respectively the complete graph, the complete bipartite,
the path and the cycle on $n$ vertices. For graph notation and terminology undefined here, we refer
the reader to \cite{Bondy-Murty2008}.

\section{Extensions of two classic inequalities}
\label{sec:classic-inequalities}
\subsection{An extension of Hong's inequality}

In 1988, Hong \cite{H88} proved that $\lambda(G)\leq \sqrt{2m-n+1}$ if $G$ is connected.
In fact, Hong's inequality holds for each graph without isolated vertices. This was 
emphasised in \cite{H93} later by Hong himself.

Let $s^+(G)$ ($s^-(G)$) be the sum of the squares of the positive (negative) eigenvalues of 
the adjacency matrix $A(G)$ of $G$. Elphick, Farber, Goldberg and Wocjan \cite{EFGW16} 
proposed the following conjecture:

\begin{conjecture}[Elphick-Farber-Goldberg-Wocjan \cite{EFGW16}] \label{conj:extension-Hong}
Let $G$ be a connected graph. Then $$\min\{s^+(G),s^-(G)\}\geq n-1.$$
\end{conjecture}

From the above conjecture, one can obtain an extension of Hong's inequality, i.e., $s^+(G)\leq 2m-n+1$.
In fact, the above conjecture has a stronger form that every graph with $\kappa$ components
satisfies that $\min\{s^+(G),s^-(G)\}\geq n-\kappa$.

Till now, Conjecture \ref{conj:extension-Hong} was confirmed for special classes of graphs, 
such as, bipartite graphs, regular graphs, and complete $k$-partite graphs \cite{EFGW16}, 
connected graphs with no more than $4$ vertices after blowing up (see \cite{Guo-Wang2022} for details).
A graph $G$ is said to be \emph{hyper-energetic} if the energy $\mathcal{E}(G)=\sum_{i=1}^n|\lambda_i|>2(n-1)$.
Conjecture \ref{conj:extension-Hong} was also confirmed for hyper-energetic graphs \cite[Theorem~10]{EFGW16}.
Since Nikiforov \cite{N07} proved that almost all graphs are hyper-energetic, Conjecture \ref{conj:extension-Hong}
is true for almost all graphs. For more families of graphs supporting Conjecture \ref{conj:extension-Hong}
and results on structural properties, we refer to \cite{ADDGHM23+}.

\subsection{An extension of Wilf's inequality}

Consider a graph $G$ of order $n$ with clique number $\omega$. In 1986, Wilf \cite{W86} proved that
\begin{equation} \label{eq:Wilf}
\lambda(G) \leq \Big(1 - \frac{1}{\omega}\Big) n.
\end{equation}
Strengthening Wilf's bound, Elphick and Wocjan \cite{Elphick-Wocjan2018} proposed a conjecture 
suggesting that \eqref{eq:Wilf} can be improved by substituting $\lambda(G)$ with $\sqrt{s^+(G)}$.

\begin{conjecture}[Elphick-Wocjan \cite{Elphick-Wocjan2018}] \label{conj:enhanced-Wilf-bound}
Let $G$ be a graph of order $n$ with clique number $\omega$. Then
\[
\sqrt{s^+(G)} \leq \Big(1 - \frac{1}{\omega}\Big) n.
\]
\end{conjecture}
This conjecture is exact, for example, for complete regular multipartite graphs. Based on some
numerical experiments, we suspect that equality holds in Conjecture \ref{conj:enhanced-Wilf-bound}
only when $G$ is a complete regular multipartite graph. In \cite{Elphick-Wocjan2018}, Elphick
and Wocjan proved this conjecture for various classes of graphs, including triangle-free graphs,
and for almost all graphs. They also tested against the thousands of named graphs with up to $40$
vertices from the Wolfram Mathematica database, but no counterexamples were found.
Using SageMath software we confirmed this conjecture for graphs having at most $10$ vertices.

\section{The Bollob\'as-Nikiforov Conjecture}

The other conjecture regarding to $\lambda_2$ is one due to Bollob\'as and Nikiforov \cite{BN07}, 
which can be seen as a spectral Tur\'an-type conjecture. The starting point of extremal graph 
theory is Mantel's theorem, which says that every graph on $n$ vertices contains a triangle
if $m>\left\lfloor n^2/4 \right\rfloor$. In 1970, Nosal \cite{N70} proved that $G$ contains
a triangle if $\lambda(G)>\sqrt{m}$. Since Nosal's result can imply Mantel's theorem,
it is always called the spectral Mantel Theorem.

Nosal's theorem was generalized by Nikiforov \cite{N02} to the inequality 
\[ 
\lambda(G)\leq \sqrt{2 \bigg(1 - \frac{1}{\omega(G)}\bigg) m}.
\] 
It should be mentioned that this 
result was implicitly suggested by Edwards and Elphick \cite{EE83}. The above 
inequality can imply the classic Tur\'an theorem on cliques, and also Wilf's 
inequality \eqref{eq:Wilf}. By introducing $\lambda_2(G)$ as a new parameter, 
Bollob\'as and Nikiforov \cite{BN07} proposed a more stronger spectral inequality.

\begin{conjecture}[Bollob\'as-Nikiforov \cite{BN07}]\label{Conj-BN}
Let $G$ be a graph on $m$ edges. Then
\[
\lambda^2_1(G)+\lambda^2_2(G)\leq 2\Big(1 - \frac{1}{\omega(G)} \Big) m.
\]
\end{conjecture}

Till now, Conjecture \ref{Conj-BN} was confirmed by Lin, Ning and Wu \cite{LNW21} for the 
case $\omega(G)=2$. They proved that for any triangle-free graph $G$ on $m$
edges and without isolated vertices, if $\lambda^2_1(G)+\lambda^2_2(G)=m$ then $G$ is a 
blow up of some member in $\{P_2, 2P_2, P_4, P_5\}$. Li, Sun and Yu \cite{Li-Sun-Yu2022}
generalized this result by giving an upper bound of $\lambda_1^{2k} + \lambda_2^{2k}$ 
for $\{C_{2i+1}\}_{i=1}^k$-free graphs. Additionally, by the result of Ando and 
Lin \cite{Ando-Lin2015}, we know that Conjecture \ref{Conj-BN} holds for weakly perfect 
graphs, which are graphs with equal clique number and chromatic number.

There is also another conjecture related to Conjecture \ref{Conj-BN} that strengthens it.

\begin{conjecture}[Elphick-Linz-Wocjan \cite{Elphick-Linz-Wocjan2021}] \label{conj:generalization-BN}
Let $G$ be a graph on $m$ edges. Then
\[
\sum_{i=1}^{\ell} \lambda_i^2(G) \leq 2\Big(1 - \frac{1}{\omega(G)} \Big) m,
\]
where $\ell = \min\{n^+(G), \omega (G)\}$, and $n^+(G)$ is the number {\rm (}counting multiplicities{\rm )}
of positive eigenvalues of $A(G)$.
\end{conjecture}

The conjecture was confirmed by Elphick, Linz and Wocjan \cite{Elphick-Linz-Wocjan2021} for 
weakly perfect graphs, Kneser graphs, Johnson graphs and two classes of strongly regular graphs.

It is important to note that choosing only $\ell=n^+(G)$ would result in a counterexample 
to Conjecture \ref{conj:generalization-BN}. In fact, consider the cycle $C_7$, one can 
check that $n^+(C_7) = 3$ and 
\[
\lambda_1^2(C_7) + \lambda_2^2(C_7) + \lambda_3^2(C_7) > 2\Big(1 - \frac{1}{\omega(C_7)} \Big) m.
\]

\section{Maximum spectral radius of planar (hyper)graphs}

Planar graphs have been extensively studied for a long time. Among various topics in spectral graph 
theory, the investigation of the spectral radius of planar graphs is particularly intriguing. This 
topic can be traced back at least to Schwenk and Wilson's fundamental question \cite{Schwenk-Wilson1978}:
``What can be said about the eigenvalues of a planar graph?"

In 1988, Hong \cite{H88} established the first significant result that for a planar graph $G$, 
the spectral radius satisfies the inequality $\lambda(G)\leq\sqrt{5n-11}$, using Hong's inequality as 
mentioned in Section \ref{sec:classic-inequalities}. Subsequently, Cao and Vince \cite{Cao-Vince1993} 
improved Hong's bound to $4+\sqrt{3n-9}$, while Hong \cite{Hong1995} himself improved it further 
to $2\sqrt{2} + \sqrt{3n - 15/2}$, and Ellingham and Zha \cite{Ellingham2000} to $2 + \sqrt{2n - 6}$. 
Additionally, Boots and Royle \cite{Boots-Royle1991} and independently, Cao and Vince \cite{Cao-Vince1993}, 
conjectured that $P_2\vee P_{n-2}$ attains the maximum spectral radius among all planar graphs 
on $n \geq 9$ vertices. Recently, significant progress has been made on the conjecture. Tait 
and Tobin \cite{Tait-Tobin2017} confirmed the conjecture for sufficiently large $n$. It is 
worth noting that Guiduli announced a proof of the conjecture for large $n$ in his Ph.D. 
Thesis (see \cite{G96} and comments from \cite{Tait-Tobin2017}). However, the conjecture 
remains open for small values of $n$.

\begin{conjecture}[Boots-Royle 1991 \cite{Boots-Royle1991} and independently Cao-Vince 1993 \cite{Cao-Vince1993}] \label{conj:planar-graphs}
Among all planar graph on $n\geq 9$ vertces, $K_2\vee P_{n-2}$ attains the maximum spectral radius.
\end{conjecture}

Extending the investigations of graph case, Ellingham, Lu and Wang \cite{ELW22} studied a hypergraph
analog of the Cvetkovi\'c--Rowlinson conjecture which states that among all outerplanar graphs on
$n$ vertices, $K_1 \vee P_{n-1}$ attains the maximum spectral radius. Given a hypergraph $H$,
the \emph{shadow} of $H$, denoted by $\partial (H)$, is the graph $G$ with $V (G) = V (H)$ and
$E(G) = \{uv : uv\in e\ \text{for some}\ e\in E(H)\}$. For a hypergraph $H$, if each edge of $H$
contains precisely $r$ vertices, then $H$ is called \emph{$r$-uniform}. The \emph{spectral radius} 
of an $r$-uniform hypergraph $H$ is defined as
\[
r! \max_{ \bm{x}\in\mathbb{R}^n,\,\|\bm{x}\|_r = 1 } \sum_{\{i_1,\ldots,i_r\}\in E(H)} x_{i_1} x_{i_2} \cdots x_{i_r},
\]
where $\mathbb{R}^n$ is the set of real vectors of dimension $n$ and $\|\bm{x}\|_r$ is the 
$\ell^r$-norm of $\bm{x}$.

A $3$-uniform hypergraph $H$ is called \emph{outerplanar} (\emph{planar}) if its shadow has an
outerplanar (planar) embedding such that each edge of $H$ is the vertex set of an interior
triangular face of the shadow. In \cite{ELW22}, Ellingham, Lu and Wang proved that for sufficiently
large $n$, the $n$-vertex outerplanar $3$-uniform hypergraph of maximum spectral radius is the
unique $3$-uniform hypergraph whose shadow is $K_1 \vee P_{n-1}$. In particular, they proposed 
a conjecture that serves as a hypergraph counterpart to Conjecture \ref{conj:planar-graphs}.

\begin{conjecture}[Ellingham-Lu-Wang \cite{ELW22}]
For large enough $n$, the $n$-vertex planar $3$-uniform hypergraph of maximum spectral radius
is the unique hypergraph whose shadow is $K_2\vee P_{n-2}$.
\end{conjecture}

\section{The relationship between Tur\'an theorem and spectral Tur\'an theorem}

A graph is said to have the \emph{Hereditarily Bounded Property} $P_{t,r}$ if $|E(H)|\leq t\cdot |V(H)| + r$
for any subgraph $H$ of $G$ with $|V(H)|\geq t$. In Guiduli's Ph.D. Thesis \cite{G96}, he proved a tight upper
bound on $\lambda(G)$ for graphs with property $P_{t,r}$.

\begin{theorem}[\cite{G96}]
Let $t\in\mathbb{N}$ and $r\geq -\binom{t+1}{2}$. If $G$ is a graph on $n$ vertices with property 
$P_{t,r}$, then
\[
\lambda (G) \leq \sqrt{tn} + \sqrt{t(t+1) + 2r} + \frac{t-1}{2},
\]
and asymptotically,
\[
\lambda(G) \leq \sqrt{tn} + \frac{t-1}{2} + o(1).
\]
Furthermore, the asymptotic bound is tight.
\end{theorem}

A natural question is to ask for a generalization of this property $P_{t,r}$, where the hereditary 
bound on the number of edges is not linear.

\begin{problem}[Guiduli \cite{G96}]
If for a graph $G$, $|E(H)|\leq c\cdot |V(H)|^2$ holds for all subgraph $H$ of $G$,
does it follows that $\lambda(G)\leq 2c\cdot |V(G)|$\,{\rm ?} What can be said if the
exponent of $2$ is replaced by some other constant less than $2$\,{\rm ?}
\end{problem}

It is worth noting that $\lambda (G) \leq \sqrt{2c}\cdot |V(G)|$ is a trivial bound from the well-known 
inequality $\lambda (G) \leq \sqrt{2 |E(G)|}$. Moreover, the constant $2c$ would be best possible,
as seen by Wilf's inequality. If this problem were true, then the spectral Erd\H os-Stone-Simonovits
theorem \cite{G96,Nikiforov2009-2} would be a consequence of the Erd\H os-Stone-Simonovits theorem \cite{Erdos1946,Erdos1966},
and Wilf's inequality \eqref{eq:Wilf} would follow from Tur\'an's theorem \cite{Turan1961}.

\section{Tight spectral conditions for a cycle of given length}

In spectral graph theory, it is very natural to ask the following problem: Determine tight spectral radius
conditions for the existence of a cycle of length $\ell$ in a graph of order $n$ for $\ell \in[3,n]$.
This problem has two aspects, i.e., the case of short cycles and the case of long cycles.

For any given graph $H$, denote by $\spex (n,H)$ the maximum spectral radius of $n$-vertex graph $G$
containing no the subgraph $H$, and by $\SPEX (n,H)$ the class of extremal graphs $G$ such that
$\lambda(G) = \spex (n,H)$. For example, when $\ell=3$, from Nosal's theorem, one can obtain $\spex (n,C_3) = \sqrt{\lfloor n^2/4 \rfloor}$.
When $n$ is odd, Nikiforov \cite{Nikiforov2007} showed $\SPEX (n,C_4) = \{K_1\vee (\frac{n-1}{2}K_2)\}$;
For the case $n$ is even, confirming a conjecture in \cite{Nikiforov2009}, Zhai and Wang \cite{Zhai-Wang2012}
showed that $\SPEX (n,C_4) = \{K_1\vee (K_1\cup \frac{n-2}{2}K_2)\}$. In 2010, Nikiforov \cite{Nikiforov2010}
proposed the following conjecture: Every graph on sufficiently large order $n$ contains a $C_{2k+2}$
if $\lambda(G)\geq \lambda(S^+_{n,k})$, unless $G=S^+_{n,k}$ where $k\geq 2$ and $S^+_{n,k}$ is
obtained from $S_{n,k}$ by adding an edge. Zhai and Lin \cite{ZL20} confirmed this conjecture
for $k=2$. Very recently, Cioab\u{a}, Desai and Tait \cite{CDT22+} announced a 
complete proof of Nikiforov's conjecture. However, the problem of determining tight spectral 
conditions for cycles of given length $\ell \in[3,n]$ is still wide open. A refined version 
of Nikiforov's conjecture can be found in \cite{LN23-1}.

\begin{problem}[A refined version of Nikiforov's even cycle conjecture \cite{LN23-1}]
For any integer $k\geq 3$, determine the infimum $\alpha:=\alpha(k)$ such that every graph of 
order $n=\Omega(k^{\alpha})$ {\rm (}where $\Omega(k^{\alpha})$ means there exists some 
constant $c$ which is not related to $k$ and $n$, such that $n\geq ck^{\alpha}${\rm )}
satisfying $\lambda(G)>\lambda(S^+_{n,k})$ contains a $C_{2k+2}$.
\end{problem}

Recently, Zhai, Lin, and Shu \cite{Zhai-Lin-Shu2021} investigated the existence of short 
consecutive cycles in fixed-size graphs and put forward the following conjecture.

\begin{conjecture}[Zhai-Lin-Shu \cite{Zhai-Lin-Shu2021}]
Let $k$ be a fixed positive integer and $G$ be a graph of sufficiently large size $m$ 
without isolated vertices. If
\[
\lambda (G) \geq \frac{k-1+\sqrt{4m-k^2+1}}{2},
\]
then $G$ contains a cycle of length $t$ for every $t\leq 2k+2$, unless $G\cong S_{m/k + (k+1)/2, k}$.
\end{conjecture}

For the case $k=2$, the conjecture has been confirmed, see \cite{Zhai-Lin-Shu2021}, \cite{Min-Lou-Huang2021} 
and \cite{Sun-Li-Wei2023} for further details.

\section{Nikiforov's problem on consecutive cycles}

A \emph{Hamilton cycle} in a graph $G$ is a cycle passing through all the vertices of $G$. If it exists,
then $G$ is called \emph{Hamiltonian}. Maybe the most famous theorem in Hamiltonian graph theory is
Dirac's theorem \cite{D52}, which states that every graph on $n\geq 3$ vertices has a Hamilton cycle
if every vertex has degree at least $n/2$. In 1971, Bondy \cite{B71} introduced the concept of pancyclicity
of graphs. Let $G$ be a graph on $n$ vertices. We say that $G$ is pancyclic, if $G$ contains each cycle 
of length $\ell$, where $\ell\in [3,n]$. Extending Ore's condition \cite{O60}, Bondy \cite{B71} proved 
that every Hamiltonian graph on $n$ vertices is pancyclic if $e(G)\geq n^2/4$, unless $G$ is isomorphic 
to $K_{n/2, n/2}$ where $n$ is even. If one drops the condition that $G$ is Hamiltonian in Bondy's 
theorem, the phenomenon of consecutive cycles still persists, i.e, there are all cycles of length 
$\ell \in [3, \lfloor (n+3)/2 \rfloor]$ in a graph $G$ if $e(G)\geq n^2/4$. This theorem may be due 
to Woodall and independently, due to Kopylov (see also Bollob\'as \cite[Corollary~5.4]{B78}).

In 2008, Nikiforov \cite{N08} considered a spectral analog of the above theorem.

\begin{problem}[Nikiforov \cite{N08}] \label{prob:consecutive-cycles}
What is the maximum $C$ such that for all positive $\varepsilon<C$ and sufficiently large $n$, every
graph $G$ of order $n$ with $\lambda (G) > \sqrt{\lfloor n^2/4 \rfloor}$ contains a cycle of length
$\ell$ for every integer $\ell\leq (C-\varepsilon)n$.
\end{problem}

Nikiforov \cite{N08} firstly showed $C\geq 1/320$ by the method of successively deleting
the least component of eigenvector with respect to spectral radius of a graph. Ning and
Peng \cite{NP20} improved it to $C\geq 1/160$. Later, Zhai and Lin \cite{ZL23} proved some
spectral result for theta graphs, and a direct main corollary is that $C\geq 1/7$. At the
same time, Li and Ning \cite{LN23+} showed that $C\geq 1/4$, by using some ideas from Ramsey
Theory \cite{ALPZ20+} and theorems on parity of cycles in graphs \cite{VZ77}. Li and Ning's
result was immediately used in \cite{ZZ23} to attack another problem in spectra graph theory.
On the other hand, Nikiforov \cite{N08} constructed the class of graphs
$G=K_s\vee (n-s)K_1$ where $s = (3-\sqrt{5})n/4$ from which one can find $C\leq (3-\sqrt{5})/2$.
Till now, Problem \ref{prob:consecutive-cycles} is still open.

\section{Graph toughness from Laplacian eigenvalues}

Let $c(G)$ denote the number of components of a graph $G$. The \emph{toughness} $t(G)$ of $G$ is 
defined by
\[
t(G) := \min\left\{\frac{|S|}{c(G-s)}\right\},
\]
in which the minimum is taken over all proper subsets $S \subset V(G)$ such that $c(G - S) > 1$. 
A graph $G$ is called $t$-tough if $t(G) \geq t$.

In 1973, Chv\'atal \cite{Chvatal1973} introduced the concept of graph toughness, which has close 
connections to a variety of graph properties such as connectivity, Hamiltonicity, pancyclicity, 
factors, and spanning trees (see \cite{Bauer-Broersma2006}). The study of toughness from 
eigenvalues was initiated by Alon \cite{Alon1995}, who showed that for any connected $d$-regular graph $G$,
\[
t(G) > \frac{1}{3} \Big( \frac{d^2}{(d + \lambda')\lambda'} - 1 \Big),
\]
where $\lambda'$ is the second largest absolute eigenvalue. Brouwer \cite{Brouwer1995} 
discovered a better bound and showed that $t(G) > d/\lambda' - 2$ for a connected $d$-regular 
graph $G$. He mentioned in \cite{Brouwer1995} that the bound might be able to be improved 
to $d/\lambda' - 1$ and then proposed the exact conjecture as an open problem 
in \cite{Brouwer1995,Brouwer1996}. In 2021, the conjecture has been proved by Gu \cite{Gu2021}.

Recently, Haemers \cite{Haemers2020} proposed studying lower bounds on $t(G)$ in terms
of the eigenvalues of the Laplacian matrix $L(G)$. He also made the following conjecture.

\begin{conjecture}[Haemers \cite{Haemers2020}]\label{conj:tougu-form-Laplacian}
Let $G$ be a connected graph on $n$ vertex with minimum degree $\delta$. Then
\[
t(G) \geq\frac{\mu_{n-1}}{\mu_1 - \delta}.
\]
\end{conjecture}

For a connected $d$-regular graph $G$, this conjecture implies that $t(G)\geq \frac{d-\lambda_2}{-\lambda_n}$,
which is stronger than Brouwer's toughness conjecture. The bound in Conjecture \ref{conj:tougu-form-Laplacian}
is tight in case $G$ is the complete multipartite graph $K_{n_1,\ldots,n_k}$ ($1 < k < n$). Indeed, assume 
$n_1 \geq n_2\geq\cdots\geq n_k$ then $t(G) = (n-n_1)/n_1$, $\mu_1 = n$ and $\mu_{n-1} = \delta = n - n_1$.

Let $S\subset V(G)$ be such that $t(G)=|S|/c(G-S)$. It was proved in \cite{Haemers2020} and \cite{Gu-Haemers2022} 
that Conjecture \ref{conj:tougu-form-Laplacian} is true in each of the following cases:
\begin{enumerate}
\item[(1)] The complement of $G$ is disconnected;
\item[(2)] All connected components of $G - S$ are singletons;
\item[(3)] The union of some components of $G - S$ has order $(n-|S|)/2$;
\item[(4)] $c(G-S) = 2$.
\end{enumerate}

In \cite{Gu-Haemers2022}, two tight lower bounds for $t(G)$ in terms of the Laplacian eigenvalues were presented,
which provided support for Conjecture \ref{conj:tougu-form-Laplacian}.

\section{Hamilton cycles in pseudo-random graphs}

Finding general conditions which ensure that a graph is Hamiltonian is a central topic in graph
theory, and researchers have devoted many efforts to obtain sufficient conditions for Hamiltonicity.

There is an old and well-known conjecture related to pseudo-random graphs in this area.
A pseudo-random graph with $n$ vertices of edge density $p$ is a graph that behaves like a truly
random graph $G(n, p)$. Spectral techniques are a convenient way to express pseudo-randomness.
An $(n,d,\lambda')$-graph is a $d$-regular graph $G$ on $n$ vertices whose second largest eigenvalue
in absolute value is at most $\lambda'$. It is known that $(n, d, \lambda')$-graphs with small
$\lambda'$ compared to $d$ possess pseudo-random properties. For more details on pseudo-random
graphs, we refer the reader to \cite{Krivelevich-Sudakov2006}. In this area, a well-known conjecture on
Hamilton cycles in an $(n,d,\lambda')$-graph can be found in \cite{Krivelevich-Sudakov2003}.

\begin{conjecture}[Krivelevich-Sudakov \cite{Krivelevich-Sudakov2003}] \label{conj:Hamilton}
There exists an absolute constant $C > 0$ such that any $(n, d, \lambda')$-graph with $d/\lambda' \geq C$
contains a Hamilton cycle.
\end{conjecture}

In \cite{Krivelevich-Sudakov2003}, Krivelevich and Sudakov proved a result that $(n,d,\lambda')$-graphs 
are Hamiltonian, provided
\[
\frac{d}{\lambda'} \geq \frac{1000\cdot \log n(\log\log\log n)}{(\log\log n)^2}
\]
for sufficiently large $n$. In recent work by Glock, Correia and Sudakov \cite{Glock-Correia-Sudakov2023}, 
progress has been made towards Conjecture \ref{conj:Hamilton} in two significant ways. Firstly, they 
improved the Krivelevich and Sudakov's bound above by showing that for some constant $C>0$, 
$d/\lambda' \geq C\cdot (\log n)^{1/3}$ already guarantees Hamiltonicity. Secondly, they established 
that for any constant $\alpha>0$, there exists a constant $C>0$ such that any $(n,d,\lambda')$-graph 
with $d\geq n^\alpha$ and $d/\lambda' \geq C$ contains a Hamilton cycle.

Let us remark that there exist three additional conjectures that are related to Conjecture \ref{conj:Hamilton}. 
The first one has to do with the concept of $f$-connected which is a generalization of the traditional notion 
of connectedness. For a graph $G$, a pair $(A,B)$ of proper subsets of $V(G)$ is called a \emph{separation} 
of $G$ if $A\cup B = V(G)$ and $G$ has no edge between $A\setminus B$ and $B\setminus A$. Let 
$f: \mathbb{N}\setminus \{0\} \to \mathbb{R}$ be a function, $G$ is called \emph{$f$-connected} if every 
separation $(A, B)$ of $G$ with $|A\setminus B| \leq |B\setminus A|$ satisfies $|A\cap B| \geq f(|A\setminus B|)$.
In 2006, Brandt, Broersma, Diestel and Kriesell \cite{Brandt-Broersma2006} conjectured that there exists a 
function $f(k) = O(k)$ (here $f(k) = O(k)$ means there exists some absolute constant $c>0$ such that 
$f(k) < ck$ for large $k$) such that every $f$-connected graph of order $n\geq 3$ is Hamiltonian. 
As asserted in \cite{Brandt-Broersma2006}, if this conjecture were true, it would imply Conjecture \ref{conj:Hamilton}.

The second one is related to Chv\'atal's toughness conjecture \cite{Chvatal1973}. In 1973, Chv\'atal
\cite{Chvatal1973} conjectured that there is a constant $t$ such that every $t$-tough graph is Hamiltonian. 
In \cite{Brandt-Broersma2006}, the authors demonstrated that the Brandt-Broersma-Diestel-Kriesell's conjecture
\cite{Brandt-Broersma2006} can be derived from Chv\'atal's toughness conjecture. Therefore, if Chv\'atal's 
conjecture were proven to be true, it would imply the validity of Conjecture \ref{conj:Hamilton}.

The third conjecture relates to the Laplacian eigenvalues of graphs. Gu (c.f.\,\cite{Gu-Haemers2022}) 
conjectured that there exists a positive constant $C < 1$ such that if $\mu_{n-1}/\mu_1 \geq C$ 
and $n \geq 3$, then $G$ contains a Hamilton cycle. It is evident that for a $d$-regular graph 
$G$, we have $\mu_{n-i+1}(G) = d - \lambda_i(G)$ for $i=1,2,\ldots,n$. Moreover, it can be 
observed that Gu's conjecture above also implies Conjecture \ref{conj:Hamilton}.

\section{A spectral problem on counting subgraphs}

Mantel's theorem, a well-known result, states that a graph with $n$ vertices and $\lfloor n^2/4 \rfloor + 1$ 
edges contains a triangle. Strengthening Mantel's theorem, Rademacher (c.f.\,\cite{Erdos1955}) proved 
that there are at least $\lfloor n/2 \rfloor$ copies of a triangle. Erd\H os \cite{Erdos1962-1,Erdos1962-2} 
further generalized the result by proving that if $q < cn$ for some small constant $c$, then 
$\lfloor n^2/4 \rfloor + q$ edges guarantees at least $q\lfloor n/2 \rfloor$ triangles. He also 
conjectured that the same result holds for $q < n/2$, which was later proved by Lov\'asz and 
Simonovits \cite{Lovasz1983}. 

In 2010, Mubayi \cite{Mubayi2010} extended these theorems to color-critical graphs, 
which are graphs whose chromatic number can be decreased by removing some edges.
Let $T_{n,k}$ denote the Tur\'an graph on $n$ vertices, which is the complete $k$-partite 
graph with parts of size $\lfloor n/k \rfloor$ or $\lceil n/k \rceil$.

\begin{theorem}[\cite{Mubayi2010}]
Let $k \geq 2$ and $F$ be a color-critical graph with chromatic number $\chi (F) = k+1$. There
exists $\delta = \delta_F > 0$ such that if $n$ is sufficiently large and $1 \leq q < \delta n$,
then every $n$-vertex graph with $e(T_{n,k}) + q$ edges contains at least $q\cdot c(n,F)$ copies of $F$,
where $c(n,F)$ is the minimum number of copies of $F$ in the graph obtained from $T_{n,k}$ by adding one edge.
\end{theorem}

Motivated by Mubayi's result above, Ning and Zhai \cite{NZ23} proposed to study 
the spectral analog of Mubayi's theorem.

\begin{problem}[Ning-Zhai \cite{NZ23}] \label{prob:supersaturation-problem}
{\rm (i)} {\rm (}The general case{\rm )} Find a spectral version of Mubayi's result.

{\rm (ii)} {\rm (}The critical case{\rm )} For $q=1$ {\rm (}where $q$ is defined as in Mubayi's theorem{\rm )},
find the tight spectral versions of Mubayi's result when $F$ is some particular color-critical subgraph,
such as triangle, clique, book, odd cycle or even wheel, etc.
\end{problem}

In \cite{NZ23}, Ning and Zhai studied the fundamental cases of triangles for Problem \ref{prob:supersaturation-problem}.
In \cite {Ning-Zhai2022}, they also studied some special bipartite case, i.e., the
quadrilaterals case.

\section{Extreme eigenvalues of nonregular graphs}

Regular graphs are a well-studied class of graphs, but for nonregular graphs where not 
all vertices have equal degrees, it is possible to quantify how close they 
are to regularity using various measures. One such measure is the difference between 
the maximum degree $\Delta(G)$ and the largest eigenvalue $\lambda(G)$ of a graph $G$. 
It is a well-known fact that $\lambda(G)\leq\Delta(G)$ for any connected graph $G$, 
and equality holds if and only if $G$ is regular. Therefore, we can use the difference 
$\Delta(G)-\lambda(G)$ as a measure of irregularity of a nonregular graph $G$. It 
is natural to ask how small this difference can be for nonregular graphs. This 
question has attracted the interest of many researchers in the past few 
decades \cite{Cioaba2007-1,Cioaba2007-2,Nikiforov2018,Shi2007,Shi2009,Stevanovic2004,Zhang2005,Zhang2021}.

Let $\mathcal{G}(n,\Delta)$ denote the set of graphs attaining the maximum spectral radius
among all connected nonregular graphs with $n$ vertices and maximum degree $\Delta$, and let
$\lambda(n, \Delta)$ denote the maximum spectral radius. For a graph $G\in\mathcal{G}(n,\Delta)$,
Liu, Shen and Wang \cite{Liu-Shen-Wang2007} investigated the order of magnitude of
$\Delta-\lambda(G)$, and posed the following conjecture: For each fixed $\Delta$ and
$G\in\mathcal{G}(n,\Delta)$, the limit of $n^2(\Delta-\lambda(G))/(\Delta-1)$ exists. Furthermore,
\[
\lim_{n\to\infty} \frac{n^2(\Delta-\lambda(G))}{\Delta-1} = \pi^2.
\]

This conjecture is trivially true for $\Delta = 2$, and the condition that $\Delta$ is fixed
is crucial (see \cite{Liu2022} for details). Recently, the first author of this paper gave a 
negative answer to the above conjecture by showing that the limit superior is at most $\pi^2/2$ (see \cite{Liu2022}).

Although Liu-Shen-Wang's conjecture is not true, we can still ask what is the exact leading 
term of $\Delta-\lambda(n,\Delta)$. Based on some numerical experiments and heuristic arguments, 
the following conjecture was presented.

\begin{conjecture}[Liu \cite{Liu2022}]
Let $\Delta\geq 3$ be a fixed integer and $G\in\mathcal{G}(n,\Delta)$. Then
the limit of $n^2 (\Delta-\lambda(G))$ always exists. Furthermore,

\begin{enumerate}
\item[$(1)$] if $\Delta$ is odd, then
\[
\lim_{n\to\infty} \frac{n^2(\Delta-\lambda(G))}{\Delta-1} = \frac{\pi^2}{4}.
\]

\item[$(2)$] if $\Delta$ is even, then
\[
\lim_{n\to\infty} \frac{n^2 (\Delta-\lambda(G))}{\Delta-2} = \frac{\pi^2}{2}.
\]
\end{enumerate}
\end{conjecture}

For $\Delta=3$ and $\Delta=4$, the conjecture was confirmed by Liu \cite{Liu2022}. However,
for the general $\Delta$, it seems to be difficult to solve it.

On the other hand, it is intuitive that the graphs attaining the maximum spectral radius among all
connected nonregular graphs with prescribed maximum degree must be close to regular graphs. In particular, Liu
and Li \cite{Liu-Li2008} posed the following conjecture: Let $3\leq\Delta\leq n-2$ and $G\in\mathcal{G}(n,\Delta)$.
Then $G$ has degree sequence $(\Delta,\ldots,\Delta,\delta)$, where
\[
\delta =
\begin{cases}
\Delta-1, & n\Delta\ \text{is odd}, \\
\Delta-2, & n\Delta\ \text{is even}.
\end{cases}
\]

In spite of the statement of the conjecture above seems intuitive, it is challenging to 
either prove or disprove it, even for small values of $\Delta$. One reason for the 
difficulty of this conjecture, as well as Liu-Shen-Wang's conjecture, is that graphs 
with bounded degree are sparse graphs whose spectral radius is bounded by a constant. 
Therefore, many tools from spectral graph theory are not effective.

By analyzing the structural properties of the extremal graphs, Liu \cite{Liu2022} confirmed Liu-Li's 
conjecture above for small $\Delta$. However, we cannot expect an affirmative answer to Liu-Li's conjecture 
for general $\Delta$. In fact, there is some evidence to support the following speculation.

\begin{conjecture}[Liu \cite{Liu2022}] \label{conj:new-degree-sequence}
Let $G\in\mathcal{G}(n,\Delta)$. For each fixed $\Delta$ and sufficiently large $n$, 
$G$ has degree sequence $(\Delta,\ldots,\Delta,\delta)$, where
\[
\delta =
\begin{cases}
\Delta-1, & \Delta\ \text{is odd},\ n\ \text{is odd}, \\
1, & \Delta\ \text{is odd},\ n\ \text{is even}, \\
\Delta-2, & \Delta\ \text{is even}.
\end{cases}
\]
\end{conjecture}

Liu confirmed Conjecture \ref{conj:new-degree-sequence} for $\Delta = 3$ and 
$\Delta = 4$ in \cite{Liu2022}.

\section{Graphs with a prescribed average degree}

While a significant amount of research has focused on finding the maximum spectral radius 
of graphs under specified conditions, there has been relatively less work done on 
determining the minimum spectral radius. 

For given $n$ and $m$, let $\mathcal{H}(n,m)$ denote the set of connected graphs on $n$ 
vertices and $m$ edges. In this section, we talk about the graphs in $\mathcal{H}(n,m)$ with 
minimum spectral radius. It is a challenging task to identify the exact graph. Nevertheless, 
the following two questions which, if true, would provide significant structural insights into 
the extremal graphs.

In 1993, Hong \cite{H93} posed a problem concerning the maximum degree and minimum degree 
of extremal graphs.

\begin{problem}[Hong \cite{H93}] \label{prop:minimum-spectral-radius}
Let $G\in\mathcal{H}(n,m)$ be a connected graph minimizing $\lambda (G) - \overline{d}(G)$. 
Is it true that $\Delta (G) - \delta (G) \leq 1$\,{\rm ?}
\end{problem}

Obviously, Problem \ref{prop:minimum-spectral-radius} is true for $m=n-1$ 
(see \cite{Collatz-Sinogowitz1957,Lovasz-Pelikan1973}); for $m=n$ (the unique 
extremal graph is clearly $C_n$); for $m=n+1$ (see \cite{Simic1989}). In \cite{Cioaba2021}, 
Cioab\u{a} claimed that Problem \ref{prop:minimum-spectral-radius} 
is also true when $2m/n$ is an integer. Indeed, If $\overline{d} := 2m/n$ 
is an integer, then it is always possible to construct a $\overline{d}$-regular 
graph which clearly minimizes the quantity $\lambda (G) - \overline{d}(G)$.

In \cite{G96}, Guiduli proposed a conjecture \cite[Conjecture 5.8]{G96} which provides a different 
perspective on the extremal graphs. It is worth mentioning that the original conjecture 
of Guiduli is not true for small $m$, the following is a modified version of Guiduli's conjecture.

\begin{conjecture}[A modified version of Guiduli's conjecture]
Let $G\in\mathcal{H}(n,m)$ be a connected graph having minimum spectral radius.
Let $r$ be such that $e(T_{n, r-1}) < m \leq e(T_{n, r})$. Then there is a constant 
$\alpha >0$ such that $G$ is $r$-colorable for $m=\Omega(n^{\alpha})$.
\end{conjecture}

\section{The Bilu-Linial Conjecture}

A \emph{signed graph} $\Gamma = (G,\sigma)$ is a graph $G = (V,E)$ along with a function $\sigma: E \to \{+1, -1\}$ 
that assigns a positive or negative sign to each edge. The (unsigned) graph $G$ is said to be the underlying graph
of $\Gamma$, while the function $\sigma$ is referred to as the \emph{signature} of $\Gamma$. 
The adjacency matrix $A(\Gamma)$ of a signed graph $\Gamma$ is derived from the adjacency matrix of 
its underlying graph $G$ by replacing every $1$ with $-1$ if the corresponding edge in $\Gamma$ is negative.
An important feature of signed graphs is the concept of switching equivalent. Given a signed 
graph $\Gamma = (G, \sigma)$ and a subset $U \subseteq V(G)$, let $\Gamma_U$ be the signed graph 
obtained from $\Gamma$ by reversing the signs of the edges in the edge boundary of $U$, the set 
of edges joining a vertex in $U$ to one not in $U$. The signed graph $\Gamma_U$ is said to 
be \emph{switching equivalent} to $\Gamma$. When we talk about the eigenvalues and spectral 
radius of a signed graph, we are actually referring to those of the corresponding signed adjacency matrix.

It was proved in \cite{Belardo-Cioaba-Koolen-Wang2018} that for a signed graph $\Gamma = (G, \sigma)$,
the spectral radius $\rho(G,\sigma)$ of $\Gamma$ is at most that of $G$. So, we know that, up to
switching equivalence, the signature leading to the maximal spectral radius is the all-positive
one. A natural question is to identify which signature leads to the minimum spectral radius.
This problem has important connections and consequences in the theory of expander
graphs \cite{Bilu-Linial2006}.

Bilu and Linial \cite{Bilu-Linial2006} proved that every regular graph has a signature 
with small spectral radius.

\begin{theorem}[\cite{Bilu-Linial2006}]
Every connected $d$-regular graph has a signature with spectral radius at most
$c \sqrt{d\cdot (\log d)^3}$, where $c > 0$ is some absolute constant.
\end{theorem}

Furthermore, they posed the following conjecture.

\begin{conjecture}[Bilu-Linial \cite{Bilu-Linial2006}]\label{conj:BiluLinial}
Every connected $d$-regular graph $G$ has a signature $\sigma$ with spectral
radius at most $2\sqrt{d-1}$.
\end{conjecture}

If true, this conjecture would provide a way to construct or show the existence of an 
infinite family of Ramanujan graphs, a connected $d$-regular graph with 
$\max\{|\lambda_2|, |\lambda_n|\} \leq 2\sqrt{d-1}$. In 2015, Marcus, Spielman and 
Srivastava \cite{Marcus-Spielman-Srivastava2015} made significant progress towards 
solving the Bilu-Linial conjecture.

\begin{theorem}[\cite{Marcus-Spielman-Srivastava2015}]
Let $G$ be a connected $d$-regular graph. Then there exists a signature $\sigma$ of
$G$ such that the largest eigenvalue of $\Gamma = (G, \sigma)$ is at most $2\sqrt{d-1}$.
\end{theorem}

\begin{remark}
In general, the largest eigenvalue $\lambda_1(\Gamma)$ of $\Gamma$ may not be equal to its
spectral radius because the Perron-Frobenius Theorem is valid only for nonnegative matrices.
To put it simply, the Bilu-Linial conjecture aims to limit all eigenvalues of $\Gamma = (G, \sigma)$ 
between $-2\sqrt{d-1}$ and $2\sqrt{d-1}$, while the Marcus-Spielman-Srivastava result
shows the existence of a signature where all the eigenvalues of $\Gamma = (G, \sigma)$ are
at most $2\sqrt{d-1}$.
\end{remark}

If the regularity assumption on $G$ is dropped, Gregory\footnote{The original link to Gregory's 
work is not available. Here we cite the description from \cite{Belardo-Cioaba-Koolen-Wang2018}.} 
considered the following variant of Conjecture
\ref{conj:BiluLinial}.

\begin{conjecture}[Gregory, c.f.\,\cite{Belardo-Cioaba-Koolen-Wang2018}] \label{conj:Gregory}
Let $G$ be a nontrivial graph with maximum degree $\Delta$. Then there exists a signature 
$\sigma$ such that $\rho(G,\sigma) < 2\sqrt{\Delta - 1}$.
\end{conjecture}

Furthermore, Belardo, Cioab\u{a}, Koolen and Wang \cite{Belardo-Cioaba-Koolen-Wang2018} posed the 
following question whose affirmative answer would imply Conjecture \ref{conj:Gregory}.

\begin{problem}[Belardo-Cioab\u{a}-Koolen-Wang \cite{Belardo-Cioaba-Koolen-Wang2018}]
Let $G$ be a connected graph. Is there a signature $\sigma$ such that $\rho(G,\sigma) < 2\sqrt{\rho(G) - 1}$\,\rm{?}
\end{problem}

\section{Isoperimetric problem in hypercube}

The $d$-dimensional hypercube, denoted by $Q_d$, is a $d$-regular graph on $2^d$ vertices, with each
vertex corresponding to a binary string of length $d$. The adjacency between two vertices in $Q_d$ 
occurs if and only if they differ in exactly one bit position. Thus, each vertex is connected to 
$d$ other vertices, which correspond to the vertices obtained by flipping each of its bits in turn.

In \cite{Bollobas-Lee-Letzter2018}, Bollob\'as, Lee and Letzter studied the maximum eigenvalue of
subgraphs of the hypercube $Q_d$. To be precise, they gave a partial answer to the following question
posed by Fink (c.f.\,\cite{Bollobas-Lee-Letzter2018}) and by Friedman and Tillich \cite{Friedman-Tillich2005}.

\begin{problem}[Fink, c.f.\,\cite{Bollobas-Lee-Letzter2018}, Friedman-Tillich \cite{Friedman-Tillich2005}]\label{prob:isoperimetric-problem}
What is the maximum of the largest eigenvalue of $Q_d[U]$, where $|U| = m$ and $1\leq m\leq 2^d$\,\rm{?}
\end{problem}

Bounding the maximum eigenvalue of $Q_d[U]$ is closely related to the size of the edge 
boundary of $U$. Indeed, since the largest eigenvalue of a graph is at least its average 
degree, for $Q_d[U]$, an induced subgraph of $Q_d$, we have $e(Q_d[U]) \leq \lambda(Q_d[U]) \cdot |U|/2$. 
Since $Q_d$ is $d$-regular, the edge boundary of $U$ is at least $(d - \lambda(Q_d[U])) \cdot |U|$. 
Hence, if we denote the maximum of the largest eigenvalue of $Q_d[U]$ with $|U|=m$ by $\lambda(m)$, then 
for every set of $m$ vertices of the hypercube $Q_d$, the edge boundary has size at 
least $(d - \lambda(m)) m$. In this sense, Problem \ref{prob:isoperimetric-problem} 
can be viewed as a variant of the isoperimetric problem in hypercube.

In \cite{Bollobas-Lee-Letzter2018}, several theorems were proved regarding this problem, 
but there are still many open problems.

\section{Minimum spectral radius of $K_{r+1}$-saturated graphs}

A graph $G$ is \emph{$F$-saturated} if $G$ does not contain $F$ as a subgraph  but the 
addition of any new edge to $G$ creates at least one copy of $F$. In other words, $G$ 
is $F$-saturated if and only if it is a maximal $F$-free graph. The maximum possible 
number of edges in a graph $G$ that is $F$-saturated is known as the Tur\'an number of $F$.
The study of Tur\'an numbers for various families of graphs is a cornerstone of extremal combinatorics.

On the other hand, the minimum number of edges in an $F$-saturated graph with $n$ vertices, 
denoted by sat$(n,F)$, is called the \emph{saturation number} of $F$. Saturation numbers were 
first studied by Erd\H os, Hajnal and Moon \cite{Erdos1964}. They determined the saturation number 
of $K_{r+1}$ and characterized the unique extremal graph, which is $S_{n,r-1}$. For a thorough 
account of the results known about saturation numbers, we refer the reader to a nice dynamic 
survey \cite{Currie-Faudree2021}.

Similarly to the spectral Tur\'an-type problems for clique, one can naturally ask whether the spectral
radius version of the Erd\H os-Hajnal-Moon theorem is true. In 2020, Kim, Kim, Kostochka and O \cite{Kim-Kostochka2020}
made a first progress on this problem, and posed the following conjecture.

\begin{conjecture}[Kim-Kim-Kostochka-O \cite{Kim-Kostochka2020}]
Let $G$ be a $K_{r+1}$-saturated graph on $n$ vertices. Then $\lambda(G) \geq \lambda (S_{n, r-1})$, with
equality if and only if $G \cong S_{n, r-1}$.
\end{conjecture}

For the cases $r=2$ and $r=3$, the conjecture was confirmed in \cite{Kim-Kostochka2020} and \cite{Kim-Kostochka-O-Shi-Wang2023}
respectively. Generally, it would be interesting to consider the following problem.

\begin{problem}
Given a graph $F$, what is the minimum spectral radius of an $F$-saturated graph on $n$ vertices\,{\rm ?}
\end{problem}

\section{Brouwer's Laplacian spectrum conjecture}

For a graph $G$ on $n$ vertices and $1\leq k \leq n$, let $S_k(G)$ denote the sum of the $k$ 
largest Laplacian eigenvalues of $G$, that is,
\[
S_k(G) := \sum_{i=1}^k \mu_i(G).
\]
As a variation of the Grone--Merris theorem \cite{Bai2011}, Brouwer \cite{Brouwer-Haemers2012} 
proposed the following conjecture.

\begin{conjecture}[Brouwer's conjecture \cite{Brouwer-Haemers2012}]\label{conj:Brouwer}
Let $G$ be a graph of order $n$. Then
\[
S_k(G) \leq e(G) + \binom{k+1}{2},~~k=1,2,\ldots,n.
\]
\end{conjecture}

The progress made on Brouwer's conjecture is worth mentioning. Brouwer himself used computers to verify 
the conjecture for all graphs with at most $10$ vertices \cite{Brouwer-Haemers2012}. For $k = 1$, the 
conjecture follows from the well-known inequality $\mu_1(G) \leq n$. Haemers, Mohammadian and 
Tayfeh-Rezaie \cite{Haemers-Mohammadian-Tayfeh-Rezaie2010} showed that the conjecture is true for 
all graphs when $k=2$, and recently the equality was characterized by Li and Guo \cite{Li-Guo2022}. 
The cases $k = n-1$ and $k = n$ are straightforward due to the fact that
$S_{n-1}(G) = S_n(G) = 2 e(G) \leq e(G) + \binom{n}{2}$.

Chen \cite{Chen2019} showed that if Conjecture \ref{conj:Brouwer} holds for all graphs when $k = p$,
then it holds for all graphs when $k = n-p-1$ as well, where $1 \leq p \leq (n-1)/2$. Thus,
Conjecture \ref{conj:Brouwer} also holds for all graphs when $k = n-2$ and $k = n-3$.
Rocha and Trevisan \cite{Rocha-Trevisan2014} proved that the conjecture is true for all
$k$ with $1 \leq k \leq \lfloor g/5 \rfloor$, where $g$ is the girth of the graph $G$
(the length of the smallest cycle in $G$). They also showed that the conjecture is true
for a connected graph $G$ having maximum degree $\Delta$, $p$ pendant vertices and $c$ 
cycles with $\Delta \geq c + p + 4$.

In addition, it has been proved that Brouwer's conjecture is true for several classes of
graphs (for all $k$) such as trees \cite{Haemers-Mohammadian-Tayfeh-Rezaie2010}, unicyclic
graphs \cite{Du-Zhou2012,Wang-Huang-Liu2012}, bicyclic graphs \cite{Du-Zhou2012}, threshold
graphs \cite{Haemers-Mohammadian-Tayfeh-Rezaie2010}, regular graphs \cite{Mayank2010} and
split graphs \cite{Mayank2010}. For more progress on Brouwer's Conjecture, we refer
to \cite{Blinovsky-Speranca2022,Chen2018,Cooper2021,Ganie-Pirzada2016,Ganie-Alghamdi-Pirzada2016,Ganie-Pirzada-Rather-Trevisan2020,Rocha2020}
and the references therein. However, Conjecture \ref{conj:Brouwer} remains open at large.

Recently, Li and Guo \cite{Li-Guo2022} proposed the following full Brouwer's conjecture.
Before continuing, we introduce some notation. For $1\leq k\leq n-1$, let $G_{k,r,s}$ $(r\geq 1, s\geq 0)$
be a graph of order $n=k+r+s$ consisting of a clique of size $k$ and two independent sets $\overline{K}_r$
and $\overline{K}_s$, where each vertex of $K_k$ is adjacent to all vertices in $\overline{K}_r$,
and for each vertex $v_i\in V(\overline{K}_s) = \{v_1,v_2,\ldots,v_s\}$, $N(v_i)\subsetneq V(K_k)$ ($i=1,2,\ldots,s$)
and $N(v_{i+1}) \subseteq N(v_i)$ ($i=1,2,\ldots,s-1$). Obviously, if $G_{k,r,s}$ is connected, then for
$k=1$, $G_{1,r,s}$ is the star $K_{1,n-1}$; for $k=n-1$, $G_{n-1,r,s}$ is the complete graph $K_n$.

\begin{conjecture}[Li-Guo \cite{Li-Guo2022}, The full Brouwer's conjecture] \label{conj:full-Brouwer}
Let $G$ be a graph of order $n$. Then
\[
S_k(G) \leq e(G) + \binom{k+1}{2},~~k=1,2,\ldots,n,
\]
with equality if and only if $G\cong G_{k,r,s}$ $(r\geq 1$, $s\geq 0)$.
\end{conjecture}

In \cite{Li-Guo2022}, the authors confirmed Conjecture \ref{conj:full-Brouwer} for $k\in \{1,2,n-3,n-2,n-1\}$.

\section{The Spectral Gap conjecture}

Given a graph $G$, the \emph{spectral gap} of $G$ is defined as $\lambda_1(G)-\lambda_2(G)$.
Obviously, if $G$ is connected, then $\lambda_1(G) - \lambda_2(G) > 0$.
The spectral gap is primarily investigated for the class of regular graphs, as it is established that
regular graphs with large spectral gap possess high connectivity properties. This property renders them
significant in numerous branches of theoretical computer science (see \cite[pp.\,392--394]{Cvetkovic-Doob-Sachs1995}).

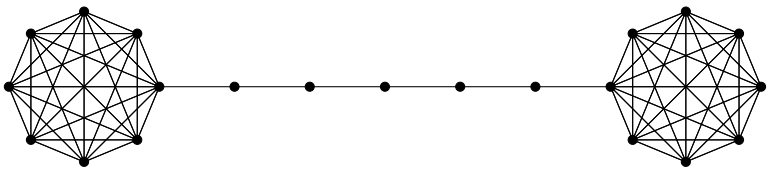
\begin{figure}[htbp]
\centering
\begin{tikzpicture}
\foreach \i in {0,1,2,...,7}
{
  \coordinate (v\i) at (45*\i:1);
  \filldraw (v\i) circle (0.06);
}
\foreach \i in {0,1,2,...,7}
    \foreach \j in {0,1,2,...,7}
       \draw (v\i) -- (v\j);
%
\foreach \i in {0,1,2,...,7}
{
  \coordinate (u\i) at ($(8,0)+(45*\i:1)$);
  \filldraw (u\i) circle (0.06);
}
\foreach \i in {0,1,2,...,7}
    \foreach \j in {0,1,2,...,7}
       \draw (u\i) -- (u\j);
%
\foreach \i in {2,3,4,5,6}
{
  \coordinate (w\i) at (\i,0);
  \filldraw (w\i) circle (0.06);
}
\draw (v0) -- (w2) -- (w3) -- (w4) -- (w5) -- (w6) -- (u4);
\end{tikzpicture}
\caption{The double kite graph $DK(8,5)$}
\label{fig:double-kite-graph}
\end{figure}

To the contrary, Stani\'c \cite{Stanic2013} suggested studying graphs with small spectral gap,
which can be viewed as an adjacency matrix version of Aldous and Fill's problem about maximizing
the relaxation time of a random walk on a connected graph (see \cite{Aksoy-Chung-Tait-Tobin2018} 
and \cite{Aldous-Fill2002} for details). In particular, Stani\'c conjectured that the minimum 
spectral gap is attained for the double kite graphs. A \emph{double kite} graph $DK(r, s)$ can 
be constructed by taking a $(s+2)$-vertex path $P_{s+2}$, two copies of the $r$-vertex complete 
graph $K_r$, and identifying one terminal vertex of $P_{s+2}$ with a vertex of one copy of $K_r$ 
and the other terminal vertex with a vertex of the other copy of $K_r$ 
(see Fig.\,\ref{fig:double-kite-graph} for an illustration).

\begin{conjecture}[Stani\'c \cite{Stanic2013}]
The spectral gap is minimized for some double kite graph over all connected graphs 
with given number of vertices.
\end{conjecture}

The conjecture has been confirmed by Stani\'c \cite{Stanic2013} for connected graphs 
with up to $10$ vertices. 

The spectral gap and the algebraic connectivity of graphs exhibit certain similarities. 
Specifically, for regular graphs, the algebraic connectivity coincides with the spectral 
gap, and connected regular graphs of degree $3$ and $4$ with minimum algebraic 
connectivity (and therefore, minimum spectral gap) are determined in \cite{Brand-Guiduli-Imrich2007} 
and \cite{Abdi-Ghorbani2023} respectively.

If we restrict our study on trees, Jovovi\'c, Koledin and Stani\'c \cite{Jovovic-Koledin-Stanic2018}
conjectured that the spectral gap is minimized for some double comet among all trees.
The \emph{double comet} $C(k, \ell)$ is a tree obtained by attaching $k$ pendant 
vertices at one end of the path $P_{\ell}$ and $k$ pendant vertices at the other 
end of the same path.

\begin{conjecture}[Jovovi\'c-Koledin-Stani\'c \cite{Jovovic-Koledin-Stanic2018}]\label{conj:spectral-grap-trees}
Among all trees of order $n$, the spectral gap is minimized for some double comet.
\end{conjecture}

The conjecture is confirmed by computer search for trees with at most $20$ 
vertices \cite{Jovovic-Koledin-Stanic2018}. There exists a unique tree that 
achieves the minimum spectral gap in all cases, and the corresponding trees 
are listed in the following table.

\begin{center}
\begin{tabular}{|c|c|c|c|c|}
\hline
$n$ & $n \leq 8$ & $9 \leq n \leq 11$ & $12 \leq n \leq 15$ & $16 \leq n \leq 20$ \\
\hline
The unique tree & $C(1, n-2) ~(\cong P_n)$ & $C(2, n-4)$ & $C(3, n-6)$ & $C(4, n-8)$ \\
\hline
\end{tabular}
\end{center}

\section{Two lower bounds on graph energy}

The graph energy is a graph-spectrum-based quantity, introduced by Ivan Gutman in the 1970s.
For a graph $G$ on $n$ vertices, the \emph{energy} $\mathcal{E}(G)$ of $G$ is defined
to be the sum of the absolute values of the eigenvalues of $A(G)$, that is,
\[
\mathcal{E}(G) := \sum_{i=1}^n |\lambda_i(G)|.
\]
This graph invariant is very closely connected to a chemical quantity known as the total $\pi$-electron
energy of conjugated hydrocarbon molecules. For results on graph energy, we refer the reader
to \cite{Li-Shi-Gutman2012}, which is a monograph summarizing the main theorems, applications
and methods regarding the adjacency energy of a graph.

The following conjecture comes from the Written on the Wall (c.f.\,\cite{Aouchiche-Hansen2010}).

\begin{conjecture}[\cite{Aouchiche-Hansen2010}] \label{conj:energy}
Let $G$ be a graph on $n$ vertices with independence number $\alpha$. Then
\[
\sum_{\lambda_i(G)>0} \lambda_i(G) \geq n - \alpha.
\]
\end{conjecture}

Note that if Conjecture \ref{conj:energy} is proven to be true, it would provide us with a concise
lower bound on the energy of graph $G$, as the left-hand side of the above inequality is precisely
equal to $\mathcal{E}(G)/2$. We utilized computer computations to verify Conjecture \ref{conj:energy}
for all graphs having at most $10$ vertices, and we did not find any counterexamples.

The well-known inertia bound \cite[Lemma 9.6.3]{Godsil-Royle2001} due to Cvetkovi\'c states that
\[
\alpha (G) \leq \min\{n - n^+, n - n^-\},
\]
where $\alpha(G)$ is the independence number of $G$, $n^+$ and $n^-$ are the numbers (counting multiplicities)
of positive and negative eigenvalues of $A(G)$ respectively. Hence, if Conjecture \ref{conj:energy} is true, 
we can promptly derive that
\begin{equation}\label{eq:energe-inertia}
\sum_{\lambda_i(G)>0} \lambda_i(G) \geq \max\{n^+, n^-\},
\end{equation}
which is also a conjecture in \cite{Aouchiche-Hansen2010}. On the ther hand, Elphick, Farber, Goldberg 
and Wocjan \cite[Lemma 5]{EFGW16} proved that
\[
s^+(G) \geq \frac{\mathcal{E}(G)^2}{4n^+} ~~ \text{and} ~~
s^-(G) \geq \frac{\mathcal{E}(G)^2}{4n^-}.
\]
Thus, Elphick pointed to us that if Conjecture \ref{conj:energy} holds true, it would imply $s^+(G) \geq n^+$ 
and $s^-(G) \geq n^-$ by \eqref{eq:energe-inertia}, which is a slightly weaker, yet analogous statement, 
compared to Conjecture \ref{conj:extension-Hong}.

The other conjecture regarding to $\mathcal{E}(G)$ is one due to Akbari and Hosseinzadeh \cite{Akbari-Hosseinzadeh2020}. 

\begin{conjecture}[Akbari-Hosseinzadeh \cite{Akbari-Hosseinzadeh2020}]\label{conj:energy-max-min-degrees}
For every graph with maximum degree $\Delta$ and minimum degree $\delta$ whose adjacency matrix is non-singular, 
$\mathcal{E}(G) \geq \Delta + \delta$ and the equality holds if and only if $G$ is a complete graph.
\end{conjecture}

In \cite{Al-Yakoob-Filipovski-Stevanovic2021}, Al-Yakoob, Filipovski and Stevanovi\'c have demonstrated 
the validity of Conjecture \ref{conj:energy-max-min-degrees} for various non-singular graphs, specifically 
those that satisfy either $n\geq \Delta+\delta$ or $|\det(A(G))|\geq\lambda(G)$ or 
$2m+n(n-1) \geq (\Delta + \delta)^2$ or $\lambda(G) - \ln\lambda(G) \geq \delta$ or 
$\Delta \leq (n-1)^{1-1/n}$.

The condition that the adjacency matrix of $G$ is non-singular is somewhat wondrous, 
Akbari and Hosseinzadeh \cite{Akbari-Hosseinzadeh2020} did not provide any motivation for the condition. 
However, as mentioned in \cite{Al-Yakoob-Filipovski-Stevanovic2021} there are obvious counterexamples 
of Conjecture \ref{conj:energy-max-min-degrees} among the graphs that have zero eigenvalues.
For example, the complete bipartite graph $K_{\Delta,\delta}$ has the adjacency spectrum:
$\sqrt{\Delta\delta}$, $-\sqrt{\Delta\delta}$ and $0$ (multiplicity $\Delta+\delta-2$). Hence, 
$\mathcal{E}(K_{\Delta,\delta}) = 2\sqrt{\Delta\delta} < \Delta + \delta$ whenever $\Delta > \delta$.

\section{Maximal $\lambda_1 + \lambda_n$ of $K_{r+1}$-free graphs}

Erd\H os put forth a conjecture that any triangle-free graph $G$ on $n$ vertices must contain a set of $n/2$ 
vertices that span at most $n^2/50$ edges, which is one of his favourite conjectures \cite{Erdos1976,Erdos1997}. 
Significant advancements have been made in various directions regarding this conjecture 
\cite{Bedenknecht2019,Erdos1994,Keevash-Sudakov2006,Norin-Yepremyan2015}. Most recently, Razborov \cite{Razborov2022} 
proved that every triangle-free graph on $n$ vertices has an induced subgraph on $n/2$ vertices with 
at most $(27/1024)n^2$ edges.

A problem with similar motivation is to determine $D_2(G)$, the minimum number of edges which has to be 
removed to make $G$ bipartite, for a triangle-free graph $G$ on $n$ vertices. A long-standing conjecture 
of Erd\H os is that at most $n^2/25$ edges need to be deleted \cite{Erdos1976}. This conjecture has been 
studied by several researchers \cite{Alon1996,Erdos-Faudree1988,Erdos-Gyori1992,Shearer1992}, and the 
most recent result was obtained by Balogh, Clemen and Lidick\'y \cite{Balogh-Clemen2021}, who proved 
that $D_2(G) \leq n^2/23.5$.

According to the Perron-Frobenius theorem, we know that $\lambda_1(G) \geq -\lambda_n(G)$. If $G$ is connected,
then equality holds if and only if $G$ is bipartite. Hence, the difference between $\lambda_1$ and $-\lambda_n$
can be viewed as a measure, how far a graph is from being bipartite. On the other hand, Brandt \cite{Brandt1998}
found a surprising connection between these two conjectures of Erd\H os and the eigenvalues of regular graphs $G$.
It was proved that
\[
\frac{\lambda_1(G) + \lambda_n(G)}{4} \cdot n \leq D_2(G).
\]

Brandt \cite{Brandt1998} also conjectured that
\[
\lambda_1(G) + \lambda_n(G) \leq \frac{4}{25} n
\]
for regular triangle-free graphs $G$ on $n$ vertices. If either of the two conjectures of Erd\H os 
were true, it would imply Brandt's conjecture (see \cite{Balogh-Clemen2022} for details).
Towards Brandt's conjecture, it was proved that $\lambda_1(G) + \lambda_n(G) \leq (3-2\sqrt{2}) n$
for regular triangle-free graphs \cite{Brandt1998}, which was shown to hold also in the non-regular 
setting by Csikv\'ari \cite{Csikvari2022}. Obviously, the quantity $\lambda_1(G) + \lambda_n(G)$ 
coincides with the smallest signless Laplacian eigenvalue of $G$ if $G$ is regular. Very recently, 
Balogh, Clemen, Lidick\'y, Norin and Volec \cite{Balogh-Clemen2022} proved that for a triangle-free 
$n$-vertex graph $G$, the smallest signless Laplacian eigenvalue of $G$ is at most $15n/94$,
which confirms Brandt's conjecture in strong form.

In general, the following problem is of interest.

\begin{problem}[Brandt \cite{Brandt1998}]
Let $r\geq 2$. How large can $\lambda_1(G) + \lambda_n(G)$ be if $G$ is a $K_{r+1}$-free graph of order $n$\,\rm{?}
\end{problem}

One can also consider a similar problem for the smallest signless Laplacian eigenvalue of graphs, see
\cite{Lima-Nikiforov-Oliveira2016} and \cite{Oboudi2022} for details.

\section{Problems on other adjacency eigenvalues of graphs}

In this section, our attention will be directed towards the adjacency eigenvalues of graphs,
excluding the largest and smallest eigenvalues.

\subsection{The third and fourth eigenvalues of graphs}

Let $G$ be a connected graph on $n$ vertices. In 1989, Powers \cite{Powers1989} presented 
an upper bound for $\lambda_i(G)$ ($1\leq i\leq n/2$) of $G$, just in terms of the order 
of $G$, i.e.,
\begin{equation}\label{eq:lambda-i}
\lambda_i(G) \leq \left\lfloor \frac{n}{i} \right\rfloor.
\end{equation}

The inequality above is notable for its simplicity and elegance. The validity of 
Inequality \eqref{eq:lambda-i} for $i\leq 2$ is clear (for $\lambda_2$, see \cite{Hong1988}), 
but unfortunately, it does not hold for $i\geq 5$ (c.f.\,\cite{Nikiforov2015}). 
Currently, the upper bound of Powers for $i\in\{3, 4\}$ remains unknown.

\begin{problem}[Powers \cite{Powers1989}] \label{prob:Powers}
Let $G$ be a graph of order $n$. Is it true that
\[
\lambda_i(G) \leq \left\lfloor \frac{n}{i} \right\rfloor
\]
for $i=3,4$\,\rm{?}
\end{problem}

Recently, Linz \cite{Linz2023} provided a counterexample to Problem \ref{prob:Powers} for 
$i=4$ by constructing a class of graphs with $\lambda_4 > n/4$. Nevertheless, it is still 
worth exploring the best possible upper bound for Problem \ref{prob:Powers}, as well as 
considering the general question posed by Hong \cite{H93}.

\begin{problem}[Hong \cite{H93}]
Find the best possible lower and upper bounds for the $i$-th eigenvalue of graphs with $n$ vertices.
\end{problem}

Let us remark that Hong's problem is related to other areas of combinatorics other than spectral 
graph theory, like the existence of symmetric Hadamard matrices, Ramsey theory and etc. 
For progress on Hong's problem we refer the reader to \cite{Nikiforov2015}.

\subsection{HL-index Conjecture}

Fowler and Pisanski \cite{Fowler-Pisanski2010-1,Fowler-Pisanski2010-2} introduced the notion of the
HL-index of a graph that is related to the HOMO--LUMO separation studied in theoretical chemistry.
This is the gap between the Highest Occupied Molecular Orbital (HOMO) and Lowest Unoccupied Molecular 
Orbital (LUMO). Let $G$ be a graph of order $n$. The \emph{HL-index} $R(G)$ of $G$ is defined as 
$R(G)=\max\{|\lambda_H(G)|, |\lambda_L(G)|\}$, where
\[
H=\Big\lfloor\frac{n+1}{2}\Big\rfloor,~~
L=\Big\lceil\frac{n+1}{2}\Big\rceil.
\]
Several bounds for this index have been proposed for some classes of graphs in 
\cite{Clemente-Cornaro2016,Fowler-Pisanski2010-1,Fowler-Pisanski2010-2,Jaklic-Fowler-Pisanski2012,Li-Li-Shi-Gutman2013}.

A graph $G$ is said to be \emph{subcubic} if its maximum degree is at most $3$. Fowler and
Pisanski \cite{Fowler-Pisanski2010-1,Fowler-Pisanski2010-2} proved that every subcubic graph
$G$ satisfies $0 \leq R(G) \leq 3$ and that if $G$ is bipartite, then $R(G) \leq \sqrt{3}$.
In 2015, Mohar \cite{Mohar2015} showed that $R(G) \leq \sqrt{2}$ for each subcubic graph $G$,
which improved the results of Fowler and Pisanski. This result is best possible since the
Heawood graph (the bipartite incidence graph of points and lines of the Fano plane) has
HL-index $\sqrt{2}$. In the same paper, Mohar also proposed the following conjecture.

\begin{conjecture}[Mohar \cite{Mohar2015}]\label{conj:HL-index}
If $G$ is a planar subcubic graph, then $R(G)\leq1$.
\end{conjecture}

The conjecture has been verified for planar bipartite graphs in \cite{Mohar2013}. Furthermore,
Mohar \cite{Mohar2016} confirmed Conjecture \ref{conj:HL-index} for all bipartite subcubic graphs
with a single exception of the Heawood graph (or a disjoint union of copies of it).

\section{Principal eigenvectors of graphs}

In the last section, we collect two conjectures on the principal eigenvectors of graphs.
For a connected graph $G$, the Perron-Frobenius theorem implies that $A(G)$ has a unique 
unit positive eigenvector corresponding to $\lambda (G)$, which is usually called 
the \emph{principal eigenvector} of $G$. 

In 2010, Cioab\u{a} \cite{Cioaba2010} presented a necessary and sufficient condition for a 
connected graph to be bipartite in terms of its principal eigenvector.

\begin{theorem}[\cite{Cioaba2010}] \label{thm:principaleigenvector}
Let $S$ be an independent set of a connected graph $G$. Then
\[
\sum_{v\in S} x_v^2 \leq \frac{1}{2}
\]
with equality if and only if $G$ is bipartite having $S$ as one color class.
\end{theorem}

Strengthening Cioab\u{a}'s result, Gregory posed the following conjecture in 2010.

\begin{conjecture}[Gregory, c.f.\,\cite{Cioaba2021}]
Let $G$ be a connected graph with chromatic number $k\geq 2$ and $S$
be an independent set. Then
\[
\sum_{v\in S} x_v^2 \leq \frac{1}{2} - \frac{k-2}{\sqrt{(k-2)^2 + 4(k-1)(n-k+1)}}.
\]
\end{conjecture}

One can check that $S_{n, k-1}$ attains equality. Let $P_r\cdot K_s$ denote the graph 
of order $(r+s-1)$ attained by identifying an end vertex of the path $P_r$ to any
vertex of the complete graph $K_s$. This graph $P_r\cdot K_s$ is called a \emph{kite graph} 
or a \emph{lollipop graph}.

The following conjecture appears in several papers \cite{Clark2022,Rucker-Gutman2002,Cioaba2021}, 
which was presented in different backgrounds.

\begin{conjecture}[R\"ucker-R\"ucker-Gutman \cite{Rucker-Gutman2002}]
Among all connected graphs on $n$ vertices, the graph $P_{n-3}\cdot K_4$ 
minimizes $\ell^1$-norm $\|\bm{x}\|_1$ of principal eigenvectors.
\end{conjecture}

The conjecture has been verified to be true for connected graphs with at most 
$10$ vertices using SageMath software.

\section*{Acknowledgment}

This paper is an invited paper of ORSC. The second author is indebted to ORSC for inviting him 
to submit a paper to the Operations Research Transactions (ORT). The authors would like to thank
Clive Elphick for his useful comments on the connection between Conjecture \ref{conj:extension-Hong} 
and Conjecture \ref{conj:energy}. The authors express sincere gratitude to Xiaofeng Gu for drawing 
attention to a related conjecture in \cite{Gu-Haemers2022} that connects to Conjecture \ref{conj:Hamilton},
and William Linz for sharing with us the reference \cite{Linz2023}. Furthermore, the authors thank
Dragan Stevanovi\'c for bringing to our attention a related conjecture in \cite{Akbari-Hosseinzadeh2020}.

\end{document}